\documentclass[12pt]{amsart}
\usepackage{amsmath, amsthm, amscd, amsfonts}

\newtheorem{theorem}{Theorem}[section]
\newtheorem{lemma}[theorem]{Lemma}
\newtheorem{proposition}[theorem]{Proposition}
\newtheorem{corollary}[theorem]{Corollary}
\theoremstyle{definition}

\theoremstyle{remark}

\numberwithin{equation}{section}
\begin{document}
\title [Automatic Continuity of $3$-Homomorphisms] {On Automatic Continuity of $3$-Homomorphisms on Banach Algebras}
\author{Janko Bra\v ci\v c}
\address{Department of Mathematics, Univ. of Ljubljana, Jadranska ul. 19, 1111 Ljubljana,
Slovenia}
\email{janko.bracic@fmf.uni-lj.si}
\author{Mohammad Sal Moslehian}
\address{Department of Mathematics, Ferdowsi University, P. O. Box 1159, Mashhad 91775, Iran;
\newline Centre of Excellence in Analysis
on Algebraic Structures (CEAAS), Ferdowsi Univ., Iran; \newline
Banach Mathematical Research Group (BMRG), Mashhad, Iran.}
\email{moslehian@ferdowsi.um.ac.ir} \subjclass[2000]{Primary 47B48,
secondary 46L05, 16Wxx.} \keywords{$3$-homomorphism, homomorphism,
continuity, Banach algebra, bounded approximate identity, second
dual, Arens regularity, $C^*$-algebra.}

\begin{abstract}
A linear map $\varphi:{\mathcal A}\to {\mathcal B} $ between
(Banach) algebras is called $3$-homomorphism if
$\varphi(abc)=\varphi(a)\varphi(b)\varphi(c)$ for each $a, b, c
\in {\mathcal A}$. We investigate $3$-homomorphisms on Banach
algebras with bounded approximate identities and establish in two ways that
every involution preserving homomorphism between $C^*$-algebras
is norm decreasing.
\end{abstract}
\maketitle

\section{Introduction.}

In \cite{HMM} the notion of $3$-homomorphism between (Banach)
algebras was introduced and some of their significant properties
were investigated. A linear mapping $\varphi:{\mathcal A}\to
{\mathcal B} $ between (Banach) algebras is called
$3$-homomorphism if $\varphi(abc)=\varphi(a)\varphi(b)\varphi(c)$
for each $a, b, c \in {\mathcal A}$. Obviously, each homomorphism
is a $3$-homomorphism. One can verify that if ${\mathcal A}$ is
unital and $\varphi$ is a $3$-homomorphism then $\psi(a) =
\varphi(1)\varphi(a)$ is a homomorphism.

There might exist $3$-homomorphisms that are not homomorphisms.
For instance, if $\psi: {\mathcal A} \to {\mathcal B}$ is a
homomorphism, then $\varphi : = -\psi$ is a
$3$-homomorphism which is not a homomorphism. As another example,
assume that ${\mathcal A}$ and ${\mathcal B} $ are Banach
algebras such that ${\mathcal A}^3=\{0\}$ and ${\mathcal B}
^3=\{0\}$. Then each linear mapping between ${\mathcal A}$ and
${\mathcal B} $ is trivially a $3$-homomorphism. One can find some non-trivial examples in \cite{HMM}.

In this paper, we investigate $3$-homomorphisms on Banach
algebras with bounded approximate identities and establish in two ways that
every involution preserving homomorphism between $C^*$-algebras is
norm decreasing.

Throughout the paper, we denote by $T^*$ the dual of a linear
mapping $T$ and by ${\mathcal A}^*$ the dual of a Banach algebra
${\mathcal A}$. We also denote the involution on a $C^*$-algebra
by $\#$.

\section{$3$-Homomorphisms on $C^*$-Algebras}

Suppose that $\varphi$ is a non-zero $3$-homomorphism from a
unital algebra ${\mathcal A}$ to ${\mathbb C}$. Then $\varphi(1)=
1$ or $-1$. Hence either $\varphi$ or $-\varphi$ is a character
on ${\mathcal A}$. If ${\mathcal A}$ is a Banach algebra, then
$\varphi$ is automatically continuous (\cite[Theorem
2.1.29]{DAL}). It may happen, however, that a $3$-homomorphism is not
continuous. For instance, let ${\mathcal A}$ be the algebra of
all $3 \times 3$ matrices having $0$ on and below the diagonal and
${\mathcal B} $ be the algebra of all ${\mathcal A}$-valued
continuous functions from $[0,1]$ into ${\mathcal A}$ with sup
norm. Then ${\mathcal B} $ is an infinite dimensional Banach
algebra in which the product of any three elements is $0$. Since
${\mathcal B} $ is infinite dimensional there are linear
discontinuous maps and so discontinuous $3$-homomorphisms from
${\mathcal B} $ into itself.

Using extreme points of the closed unit ball of a $W^*$-algebra
${\mathcal A}$ authors of \cite{HMM} showed that if ${\mathcal B} $
is a $C^*$-algebra and $\varphi:{\mathcal A}\to {\mathcal B} $ is a
weakly-norm continuous involution preserving $3$-homomorphism, then
$\Vert \varphi\Vert\leq 1$. Then they asked whether every preserving
involution $3$-homomorphism between $C^*$-algebras is continuous. We
shall establish this by applying the strategy used by L. Harris in
\cite[Proposition 3.4]{HAR} where it is shown that each
$J^*$-homomorphism between $J^*$-algebras is norm decreasing.

\begin{theorem} \label{theorem} Let $\varphi:{\mathcal A}\to {\mathcal B} $ be an involution preserving $3$-homomorphism
between unital $C^*$-algebras. Then $\|\varphi\|\leq
1$.\end{theorem}
\begin{proof} Let $a \in {\mathcal A}$, let
$\lambda>0$ and let $\lambda \notin \sigma(a^\#a)$. Then $\lambda
1_{{\mathcal A}}-a^\#a$ has an inverse $b$. Then we have the
following sequence of deductions:

\[ab(\lambda 1_{{\mathcal A}}-a^\#a)=a\]
\[(\lambda b-ba^\#a)a(\lambda b-ba^\#a)=a\]
\[\varphi(\lambda b-ba^\#a)\varphi(a)\big(\lambda\varphi(b)
-\varphi(b)\varphi(a^\#)\varphi(a)\big)=\varphi(a)\]
\[\varphi(a^\#)\,\big(\varphi(\lambda b-ba^\#a)\varphi(a)\varphi(b)\big)\,\big(\lambda
1_{{\mathcal B}}
-\varphi(a^\#)\varphi(a)\big)=\varphi(a^\#)\varphi(a)\]
\[\varphi(a^\#)\varphi(ab)\big(\lambda1_{{\mathcal B}}-\varphi(a^\#)\varphi(a)\big)=\varphi(a^\#)\varphi(a)\]
\[\lambda 1_{{\mathcal B}}-\varphi(a^\#)\varphi(ab)\big(\lambda
1_{{\mathcal B}}-\varphi(a^\#)\varphi(a)\big) =\lambda 1_{{\mathcal
B}}-\varphi(a^\#)\varphi(a)\]
\[1_{{\mathcal B}}=\frac{1}{\lambda}\big(1_{{\mathcal
B}}+\varphi(a^\#)\varphi(ab)\big)(\lambda 1_{{\mathcal B}}-
\varphi(a^\#)\varphi(a)\big)\] Since $\lambda 1_{{\mathcal
B}}-\varphi(a^\#)\varphi(a)$ is self-adjoint, it follows that
$\lambda 1_{{\mathcal B}}-\varphi(a)^\#\varphi(a)$ is invertible and
so $\lambda\notin \sigma\big(\varphi(a)^\#\varphi(a)\big)$. Hence
\[\sigma\big(\varphi(a)^\#\varphi(a)\big)\subseteq \sigma(a^\#a)\cup\{0\}.\]
Since the norm of a normal element $c$ of a $C^*$-algebra is
equal to its spectral radius $r(c)$ (see \cite[Theorem 3.2.3]{DAL}), we have
\[\|\varphi(a)\|^2=\|\varphi(a^\#a)\|=r\big(\varphi(a)^\#\varphi(a)\big)\leq
r(a^\#a) =\|a^\#a\|=\|a\|^2.\] So $\varphi$ is norm
decreasing.\end{proof}

\section{$3$-homomorphisms from algebras with bounded approximate identity}

Let ${\mathcal A}$ be a Banach algebra. Let $a\in {\mathcal A}$
and $\xi\in {\mathcal A}^*$ be arbitrary, then $a.\xi \in
{\mathcal A}^*$ and $\xi.a \in {\mathcal A}^*$ are defined by
$\langle a.\xi,x\rangle=\langle \xi,xa\rangle, x \in {\mathcal
A}$ and
$\langle \xi.a,x\rangle =\langle \xi,ax\rangle, x \in {\mathcal A}.$\\
For $\xi \in {\mathcal A}^*$ and $F \in {\mathcal A}^{**}$,
functionals $F.\xi \in {\mathcal A}^*$ and $\xi.F \in {\mathcal
A}^*$ are given by $\langle F.\xi,x\rangle =\langle
F,\xi.x\rangle, x \in {\mathcal A}$ and $\langle \xi.F,x\rangle
=\langle F,x.\xi \rangle, x \in {\mathcal A}$. Recall that the
first Arens product of $F,G \in {\mathcal A}^{**}$ is
$F\triangleleft G \in {\mathcal A}^{**}$, where
\[\langle F\triangleleft G,\xi\rangle=\langle F,G.\xi \rangle \qquad (\xi \in {\mathcal
A}^*).\]
Similarly, the second Arens product of $F,G \in {\mathcal
A}^{**}$ is $F\triangleright G \in {\mathcal A}^{**}$, where
\[\langle F\triangleright G,\xi\rangle=\langle G, \xi.F \rangle
\qquad (\xi \in {\mathcal A}^*).\] It is well known that
$({\mathcal A}^{**},\triangleleft)$ and $({\mathcal A}^{**},
\triangleright)$ are Banach algebras and that ${\mathcal A}$ is a
closed subalgebra of each of them. The algebra ${\mathcal A}$ is
Arens regular if the products $\triangleleft$ and
$\triangleright$ coincide on ${\mathcal A}^{**}$.

A left approximate identity for ${\mathcal A}$ is a net
$(e_{\alpha})\subset {\mathcal A}$ such that $
\lim_{\alpha}e_{\alpha}a=a$, for each $a \in {\mathcal A}$. If
there exists $m>0$ such that $\sup_{\alpha}\|e_{\alpha}\| \leq
m$, the left approximate identity $(e_{\alpha})$ is said to be
bounded by $m$. The definitions of a right approximate identity
and a bounded right approximate identity are similar. Recall from
\cite[Proposition 2.9.16]{DAL} that ${\mathcal A}$ has a left
approximate identity bounded by $m$ if and only if $({\mathcal
A}^{**},\triangleright)$ has a left identity of norm at most
$m$. Similarly, ${\mathcal A}$ has a right approximate identity
bounded by $m$ if and only if $({\mathcal A}^{**},\triangleleft)$
has a right identity of norm at most $m$.

\begin{lemma}
Let ${\mathcal A}$ and ${\mathcal B} $ be Banach algebras and
$\varphi: {\mathcal A} \to {\mathcal B} $ be a $3$-homomorphism.
Then, for arbitrary $a, b \in {\mathcal A}$, $\xi \in {\mathcal B}
^*$, and $F, G \in {\mathcal A}^{**}$, the following hold:

(i)\, $\varphi^*(\varphi(a)\varphi(b).\xi)=ab.\varphi^*(\xi); $

(ii)\, $\varphi^*(\xi.\varphi(a)\varphi(b))=\varphi^*(\xi).ab; $

(iii)\, $\varphi^*(\varphi(a).\xi.\varphi^{**}(F))=a.\varphi^*(\xi).F; $

(iv)\, $\varphi^*(\varphi^{**}(F).\xi.\varphi(a))=F.\varphi^*(\xi).a; $

(v)\, $\varphi^*(\xi.\varphi^{**}(F)\triangleright \varphi^{**}(G))
=\varphi^*(\xi).F\triangleright G; $

(vi)\, $ \varphi^*(\varphi^{**}(F)\triangleleft \varphi^{**}(G).\xi)
=F\triangleleft G.\varphi^*(\xi).$
\end{lemma}

\begin{proof}
(i) Let $a, b \in {\mathcal A}$, and $\xi \in {\mathcal B} ^*$ be
arbitrary. Then, for each $x \in {\mathcal A}$,
\begin{eqnarray*}
\langle \varphi^*(\varphi(a)\varphi(b).\xi),x\rangle &=&\langle
\varphi(a)\varphi(b).\xi, \varphi(x)\rangle\\
&=&\langle \xi,\varphi(x)\varphi(a)\varphi(b)\rangle \\
&=&\langle \xi,\varphi(xab)\rangle\\
&=&\langle \varphi^*(\xi), xab\rangle\\
&=& \langle ab.\varphi^*(\xi),x\rangle.
\end{eqnarray*}

(iii) Let $a \in {\mathcal A}$, $\xi \in {\mathcal B} ^*$, and $F
\in {\mathcal A}^{**}$ be arbitrary. Then, for each $x \in
{\mathcal A}$,
\begin{eqnarray*}
\langle \varphi^*(\varphi(a).\xi.\varphi^{**}(F)),x\rangle &=&
\langle \varphi(a).\xi.\varphi^{**}(F), \varphi(x)\rangle \\
&=&\langle \varphi^{**}(F),\varphi(x)\varphi(a).\xi \rangle\\
&=& \langle F,\varphi^*(\varphi(x)\varphi(a).\xi)\rangle\\
&=&\langle F,xa.\varphi^*(\xi)\rangle\\
&=&\langle a.\varphi^*(\xi).F,x\rangle.
\end{eqnarray*}

(v) Let $\xi \in {\mathcal B} ^*$, and $F, G \in {\mathcal
A}^{**}$ be arbitrary. Then, for each $x \in {\mathcal A}$,
\begin{eqnarray*}
\langle \varphi^*(\xi.\varphi^{**}(F)\triangleright
\varphi^{**}(G)),x\rangle &=& \langle
(\xi.\varphi^{**}(F)\triangleright
\varphi^{**}(G)), \varphi(x)\rangle\\
&=&\langle \varphi^{**}(G),\varphi(x).\xi.\varphi^{**}(F)\rangle\\
&=&\langle G,\varphi^*(\varphi(x).\xi.\varphi^{**}(F))\rangle\\
&=&\langle G,x.\varphi^*(\xi).F\rangle \\
&=&\langle \varphi^*(\xi).F\triangleright G,x\rangle.
\end{eqnarray*}
The proofs of (ii), (iv), and (vi) are similar.
\end{proof}

\begin{proposition}\label{prop1}
Let ${\mathcal A}$ and ${\mathcal B} $ be Banach algebras. If
$\varphi: {\mathcal A} \to {\mathcal B} $ is a $3$-homomorphism,
then $\varphi^{**}$ is a $3$-homomorphism from $({\mathcal
A}^{**},\triangleleft)$ to $({\mathcal B} ^{**},\triangleleft)$
and from $({\mathcal A}^{**},\triangleright)$ to $({\mathcal B}
^{**},\triangleright)$.
\end{proposition}

\begin{proof}
Let $F, G, H \in {\mathcal A}^{**}$ be arbitrary. Then, for each
$\xi \in {\mathcal B} ^*$,
\begin{eqnarray*}
\langle \varphi^{**}(F\triangleleft G\triangleleft H),\xi \rangle
&=&\langle F\triangleleft G\triangleleft H, \varphi^{*}(\xi)\rangle\\
&=&\langle F,G\triangleleft H.\varphi^*(\xi)\rangle\\
&=&\langle F,\varphi^*(\varphi^{**}(G)\triangleleft \varphi^{**}(H).\xi)\rangle \\
&=&\langle \varphi^{**}(F), \varphi^{**}(G)\triangleleft (\varphi^{**}(H).\xi)\rangle \\
&=&\langle \varphi^{**}(F)\triangleleft \varphi^{**}(G), \varphi^{**}(H).\xi\rangle \\
&=&\langle \varphi^{**}(F)\triangleleft \varphi^{**}(G)\triangleleft
 \varphi^{**}(H),\xi \rangle,
\end{eqnarray*}
where the second equality holds by Lemma 5.1 (vi).

The proof of the second assertion is similar.\end{proof}

\begin{theorem} Let ${\mathcal A}$ and ${\mathcal B} $ be Banach algebras and $\varphi: {\mathcal A} \to {\mathcal B} $ be a
$3$-homomorphism. If ${\mathcal A}$ has a bounded left approximate
identity, then
\[\psi(F):=\varphi^{**}(L)\triangleright \varphi^{**}(F)\;\;\;\;\;(F \in {\mathcal
A}^{**})\] defines a homomorphism from $({\mathcal A}^{**},
\triangleright)$ to $({\mathcal B} ^{**},\triangleright)$, where
$L$ is a left identity in ${\mathcal A}^{**}$. Moreover,
\[\varphi(a)=\varphi^{**}(L)\triangleright\psi(a)\;\;\;\;\;(a \in
{\mathcal A}).\] Similarly, if ${\mathcal A}$ has a bounded right
approximate identity, then
\[\psi(F):=\varphi^{**}(F)\triangleleft
\varphi^{**}(R)\;\;\;\;\;(F \in {\mathcal A}^{**})\] defines a
homomorphism from $({\mathcal A}^{**}, \triangleleft)$ to
$({\mathcal B} ^{**},\triangleleft)$, where $R$ is a left
identity in ${\mathcal A}^{**}$. Moreover,
\[\varphi(a)=\psi(a)\triangleleft \varphi^{**}(R)\;\;\;\;\;(a \in
{\mathcal A}).\]
\end{theorem}

\begin{proof} We shall prove only the first part of the
theorem. Since $\varphi$ is a $3$-homomorphism the mapping
$\varphi^{**}$ is a $3$-homomorphism from $({\mathcal A}^{**},
\triangleright)$ to $({\mathcal B} ^{**}, \triangleright)$, by
Proposition \ref{prop1}. By \cite[Proposition 2.9.16]{DAL},
there exists a left identity $L$ in $({\mathcal
A}^{**},\triangleright)$. Thus
\[\psi(F):=\varphi^{**}(L)\triangleright
\varphi^{**}(F)\;\;\;\;\;(F \in {\mathcal A}^{**})\] defines a
homomorphism from $({\mathcal A}^{**},\triangleright)$ to
$({\mathcal B} ^{**},\triangleright)$ and
\[\varphi^{**}(F)=\varphi^{**}(L)\triangleright
\psi(F)\;\;\;\;\;(F \in {\mathcal A}^{**}).\]
Since the
restriction of $\varphi^{**}$ to ${\mathcal A}$ is $\varphi$, we
have
\[\varphi(a)=\varphi^{**}(L)\triangleright
\psi(a)\;\;\;\;\;(a  \in {\mathcal A}).\]
\end{proof}

We now give another proof of Theorem \ref{theorem}:

\begin{corollary} Suppose that ${\mathcal A}$ and ${\mathcal B}$ are $C^*$-algebras. Then
each involution preserving $3$-homomorphism $\varphi: {\mathcal
A} \to {\mathcal B}$ is bounded with norm $\|\varphi\| \leq 1$.
\end{corollary}
\begin{proof} Each $C^*$-algebra has approximate identity of bound
1 (cf. \cite[Theorem 3.2.21]{DAL}). By \cite[Corollary 3.2.37]{DAL}, ${\mathcal A}$ and ${\mathcal B} $ are Arens regular
and ${\mathcal A}^{**}$ and ${\mathcal B} ^{**}$ are von Neumann
algebras equipped with the involution that are defined in the
following way (see \cite[Page 349]{DAL}). A linear involution on
${\mathcal A}^*$ is defined by
\[\langle \xi^\#,a\rangle=\overline{\langle
\xi,a^\#\rangle}\;\;\;\;(a \in {\mathcal A}, \xi \in {\mathcal
A}^*),\] and an involution on ${\mathcal A}^{**}$ is given by
\[\langle F^\#,\xi\rangle =\overline{\langle
F,\xi^\#\rangle}\;\;\;\;\;(F \in {\mathcal A}^{**}, \xi  \in
{\mathcal A}^*).\] The involution on ${\mathcal B} ^{**}$ is
defined similarly. The algebra ${\mathcal A}^{**}$ is unital,
with identity $E$  of norm $\|E\|=1$.

It is easy to check that $\varphi^{**}$ is an involution
preserving $3$-homomorphism. Namely, for each $a \in {\mathcal A}$
and each $\xi \in {\mathcal B} ^*$, we have
\[\langle \varphi^*(\xi^\#),a\rangle=\overline{\langle
\xi,\varphi(a)^\#\rangle}= \overline{\langle
\xi,\varphi(a^\#)\rangle}=\langle \varphi^*(\xi)^\#,a\rangle.\] It
follows, for each $F \in {\mathcal A}^{**}$ and each $\xi \in
{\mathcal B} ^*$, that
\[\langle
\varphi^{**}(F^\#),\xi\rangle=\overline{\langle
F,\varphi^*(\xi)^\#\rangle}= \overline{\langle
F,\varphi^*(\xi^\#)\rangle}=\langle
\varphi^{**}(F)^\#,\xi\rangle.\]

Since
\begin{eqnarray*}
\varphi^{**}(E)\varphi^{**}(F)&=&\varphi^{**}(E)\varphi^{**}(FE^2)\\
&=&\varphi^{**}(E) \varphi^{**}(F)\varphi^{**}(E)^2\\
&=&\varphi^{**}(EFE)\varphi^{**}(E)\\
&=&\varphi^{**}(F)\varphi^{**}(E),
\end{eqnarray*}
for each $F \in
{\mathcal A}^{**}$, we conclude that
$\psi(F)=\varphi^{**}(E)\varphi^{**}(F)\;\;\;(F  \in {\mathcal
A}^{**})$ defines an involution preserving homomorphism. Indeed,
for each $F \in {\mathcal A}^{**}$, we have
\[\psi(F^\#)=\varphi^{**}(E)\varphi^{**}(F^\#)=(\varphi^{**}(F)\varphi^{**}(E)^\#)^\#
=\psi(F)^\#.\] By \cite[Corollary 3.2.4]{DAL}, each involution
preserving homomorphism between $C^*$-algebras is norm-decreasing.
Thus
\[ \|\varphi^{**}(F)\|=\| \varphi^{**}(E)\psi(F)\| \leq \|
\varphi^{**}(E)\| \|F\|\;\;\;\;\; (F \in {\mathcal A}^{**}).\]
Since
\[\psi(E)=\varphi^{**}(E)^2=\varphi^{**}(E)^\#\varphi^{**}(E)\]
we have \[1\geq \|
\psi(E)\|=\|\varphi^{**}(E)^\#\varphi^{**}(E)\|=\|\varphi^{**}(E)\|^2,\]
which gives $ \| \varphi^{**}(E)\| \leq 1.$ Thus,
\[\|\varphi^{**}(F)\| \leq \| \varphi^{**}(E)\| \|F\| \leq \|F\|\;\;\;\;\;
(F \in {\mathcal A}^{**})\].
\end{proof}

{\bf Acknowledgement.}  The authors would like to thank the referees
for their valuable comments and Professor John. D. Trout for his
useful comment on Theorem 2.1.

\end{document}